\documentclass[a4paper,10pt]{article}

\usepackage{latexsym}
\usepackage{amsmath,amsthm,amsfonts}
\usepackage{graphicx}

\newcommand{\bZ}{\mathbb{Z}}

\DeclareMathOperator{\Aut}{Aut}

\theoremstyle{plain}

\newtheorem{Theorem}{Theorem}

\newtheorem{Proposition}[Theorem]{Proposition}

\newtheorem{Problem}{Problem}

\theoremstyle{remark}

\begin{document}

\title{Multiple Kronecker Covering Graphs}

\author{
  Wilfried Imrich\thanks{\texttt{Wilfried.Imrich@mu-leoben.at}},\\
  Montanuniversit\"{a}t Leoben, Austria,\\ and \\
  Toma\v{z}
  Pisanski\thanks{\texttt{Tomaz.Pisanski@fmf.uni-lj.si}},\\
  IMFM,University of Ljubljana, and University of Primorska, Slovenia
}
\date{ }

\maketitle

\begin{abstract}
A graph may be the Kronecker cover in more than one way. In this
note we explore this phenomenon. Using this approach we show that
the least common cover of two graphs need not be unique.
\end{abstract}

\section{Introduction}

A graph $\tilde G$ is said to be a covering graph over a graph $G$ if there
exists a surjective homomorphism (called a covering) $f\colon\tilde G\to G$
such that for every vertex $v$ of $\tilde G$ the set of edges incident with
$v$ is mapped bijectively onto the set of edges incident with $f(v)$.
A covering $f$ is $k$-fold if the preimage of every vertex of $G$ consists
of $k$ vertices.

To simplify the description of large graphs, the concept of \emph{voltage graphs}
and \emph{covering graphs} is generally used, see for example~\cite{GT} or~\cite{White}.

In 1982 F.T. Leighton proved in \cite{Leighton} that any two graphs
with a common universal cover have a common finite cover. It is not
hard to see that any two graphs with a common cover have a unique
maximal common cover: the universal cover. In this note we show that
the result does not extend in the opposite direction. There are
graphs with a common cover whose minimal common cover is not unique.
\section{Graphs are not determined by their Kronecker covers}

The \emph{Kronecker cover} of $G$ is the $\bZ_2$-covering graph over
$G$ with voltages $1$ on all edges (note that in this case the
direction of edges is irrelevant). We will denote the Kronecker
cover of $G$ by $KC(G)$. Alternatively, $KC(G)$ can be defined as
the tensor product of $G$ and $K_2$. See~\cite{ImrKla} for more
about graph products.

It is easy to see the following properties of $KC(G)$.
\begin{Proposition}
Kronecker covers of graphs are bipartite. If $G$ is bipartite, then $KC(G)$
consists of two copies of $G$. If $G$ is connected and non-bipartite then $KC(G)$
is connected.
\end{Proposition}
\begin{proof}
By definition, the vertex set $V(KC(G))$ is a union of two sets $V(G)$ and the
edges of $KC(G)$ connect only vertices in different copies of $V(G)$.
\end{proof}

It follows from the theory of tensor products developed by Imrich et al. that
$KC(G) = KC(G')$ does not necessarily imply $G = G'$. Recently, Imrich et al. have
determined all possibilities for a hypercube $Q_n$ to be a Kronecker cover.

\begin{figure}
\includegraphics[width=0.5\textwidth]{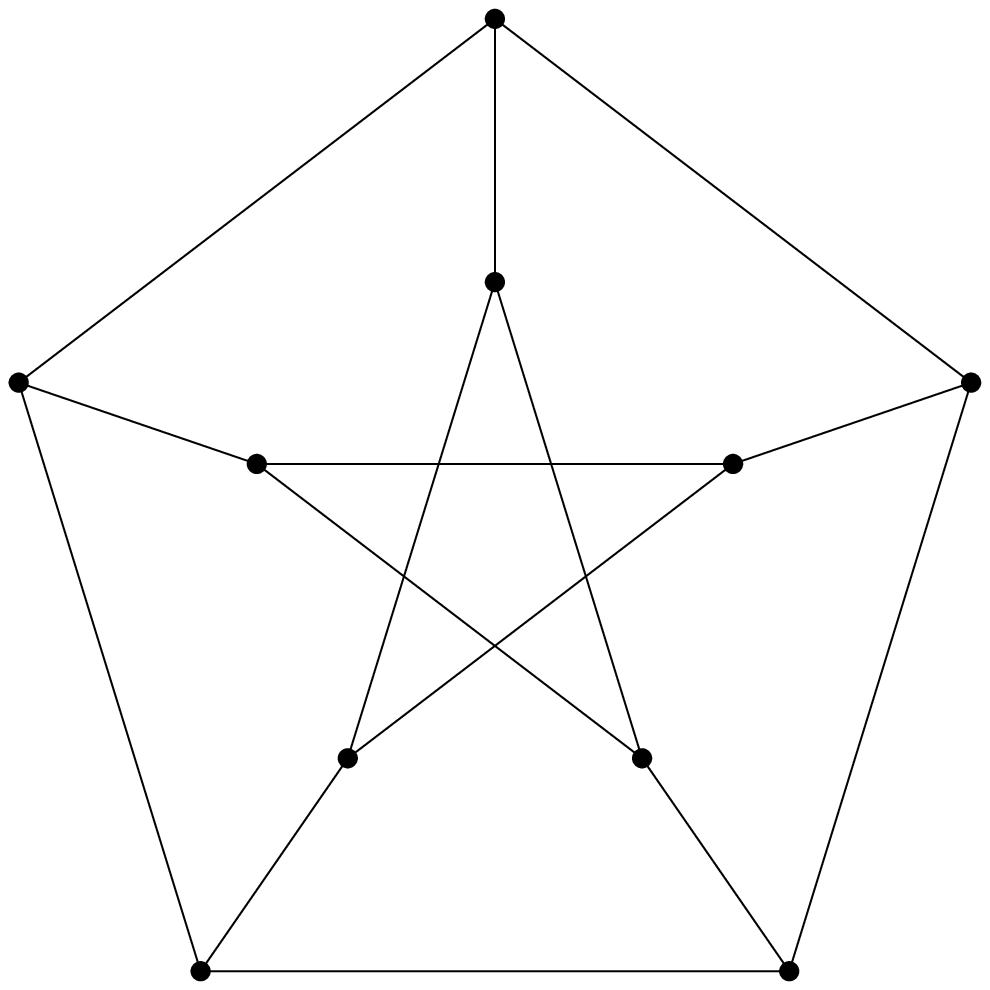}
\includegraphics[width=0.5\textwidth]{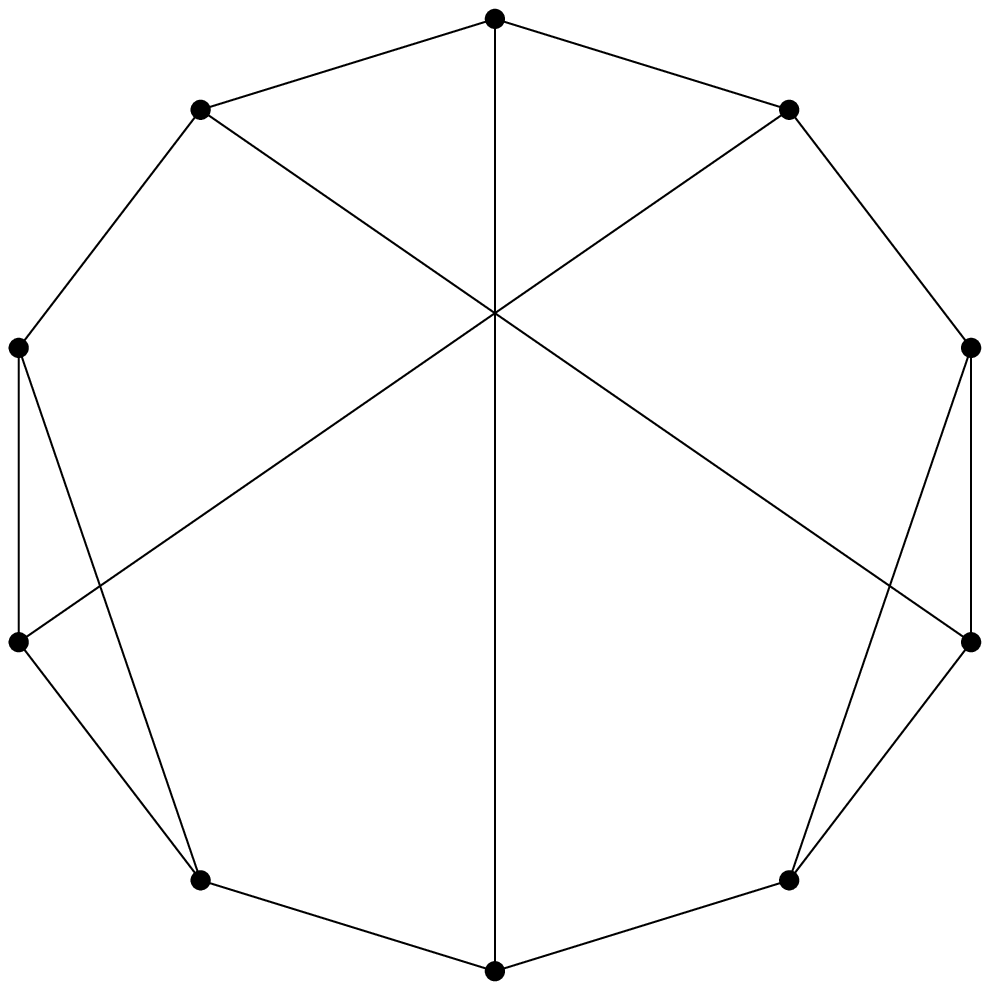}
\caption{\label{fig_1} The Desargues graph $G(10,3)$ is a Kronecker
cover of Petersen graph $G(5,2)$ and of the graph $X$.}
\end{figure}

Here we open the problem for all simple graphs.

\begin{Problem}
Given a connected, simple graph $K$, determine all simple graphs $G$ such that $K = KC(G)$.
\end{Problem}

Clearly, not all graphs can be Kronecker covers. Here is a simple criterion for
a graph $K$ to be Kronecker cover. Let $K$ be bipartite with bipartition $(V_1,V_2)$
and let $\pi \in \Aut K$ be a
fixed-point free involution such that $\pi$ interchanges the bipartition:
$\pi(V_1) = V_2.$ Furthermore, we require that for any vertex $v$ of $K$ vertices $v$ and $\pi(v)$ are
non-adjacent.
Such an automorphism is called a {\em (combinatorial) polarity}.

\begin{Proposition}
Let $K$ be a connected graph. Then $K$ is a Kronecker cover of some graph $G$ if and only
if $K$ is bipartite and there exists a polarity $\pi \in \Aut K$.
\end{Proposition}

Let $\Pi K \subset \Aut K$ denote the set of all polarities of $K$.
Clearly if $\pi$ is a polarity and $\alpha$ an arbitrary
automorphism, then $\pi^{\alpha} = \alpha \pi \alpha^{-1}$ is also a
polarity, because every automorphism either fixes or interchanges
the bipartition. Let $\pi^{\Aut K}$ denote the class $\pi^{\Aut K} =
\{\pi^{\alpha}|\alpha \in \Aut K \}$. If we define an equivalence
relation $\cong$ in $\Pi K$ so that $\pi$ is equivalent to $\pi'$ if
and only if there exists an $\alpha \in \Aut K$ such that $\pi' =
\pi^{\alpha},$ then the equivalence classes are exactly of the form
$\pi^{\Aut K}.$

\begin{Proposition}
Let $K$ be a connected graph. Then $K$ is a Kronecker cover of $k$ simple graphs if and only if
$k = |\Pi_K/{\cong}|$.
\end{Proposition}

Let us consider the case presented in Figure~\ref{fig_1}. The
Desargues graph $G(10,3)$ can be represented as a Kronecker covering
graph in two distinct ways. Let us label the vertices of Petersen
graph $G(5,2)$ in such a way that the vertices in the outer pentagon
are labeled 1,2,3,4,5 and the vertices in the inner pentagram
6,7,8,9,10 with 1 being adjacent to 6, 2 to 7, etc. The labeling of
the vertices of the bipartite $G(10,3)$ is chosen in such a way that
the vertex $i$ of Petersen lifts to a black vertex $i$ and a white
vertex $i'$. Each pair $i$ and $i'$ of vertices is antipodal in
$G(10,3).$ In order to specify the second quotient, the graph $X$ we
have to define a new polarity $\pi$ of $G(10,3)$ that tells which
black vertex $\pi(i)$ projects onto the vertex $i$ of $X$. We do
this with the aid of an involution $\alpha$ if $G(5,2)$ defined by
$\alpha = (1,8)(2,10)(3,5),(4)(6)(7)(9)$ with four fixed points by
setting $\pi(i) = \alpha(i)'$ and $\pi(i') = \alpha(i)$. If we now
identify $\pi(i')$ with $i$ we obtain a covering projection of
$G(10,3)$ onto the graph $X$ labelled in the following order along
the Hamilton cycle of the graph on the right side of Figure 1:
$\{10,1,3,4,9,6,8,5,2,7\}$ where $10-1-3$ and $8-5-2$ are the two
triangles. More generally we have the following proposition.

\begin{Proposition} Let $\alpha$ be an involution of a graph $G$ that does
not interchange the endpoints of an edge, then $\pi(i) = \alpha(i)'$
and $\pi(i') = \alpha(i)$ is a polarity of $KC(G)$.
\end{Proposition}

\begin{figure}
\includegraphics[width=0.5\textwidth]{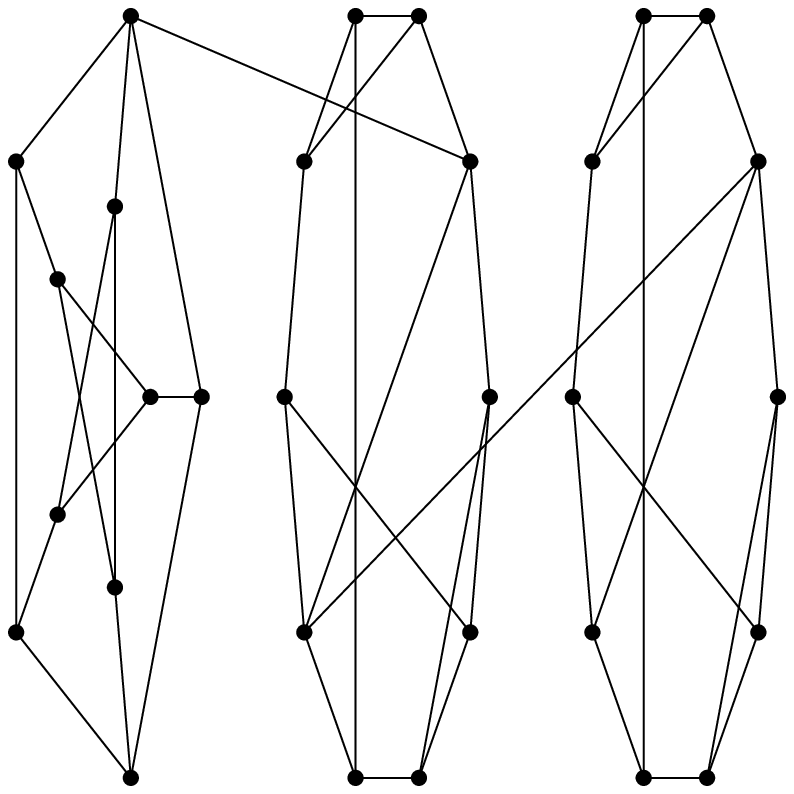}
\includegraphics[width=0.5\textwidth]{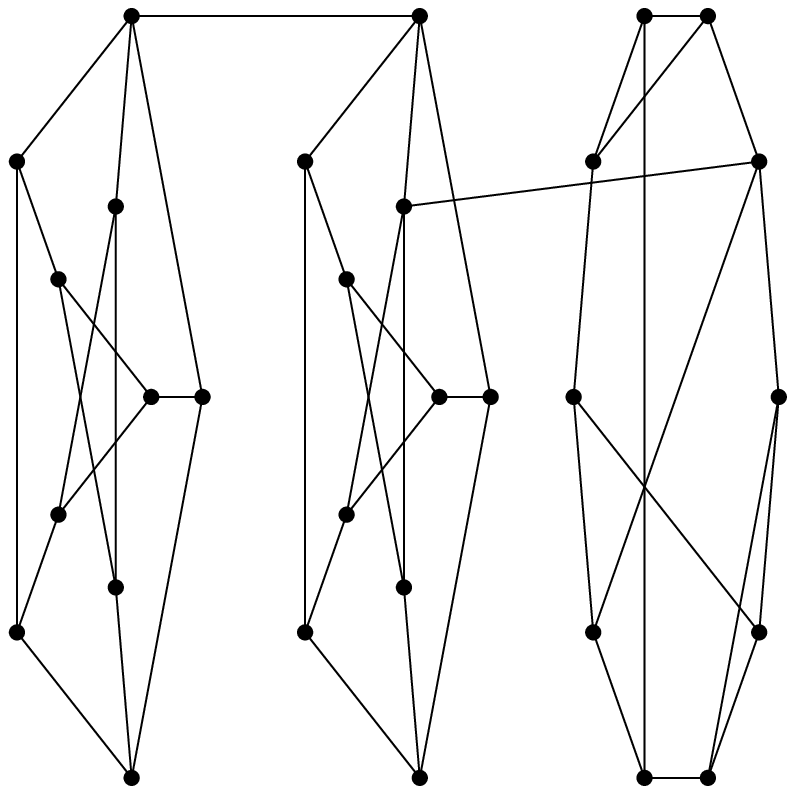}
\caption{\label{fig_2} The graphs $H_1$ and $H_2$ have more than one
common connected double cover.}
\end{figure}

\section{Least common covers are not unique}
It is well-known that if two graphs share a common cover then they
have the same universal cover, that is the largest possible
connected cover of the two graphs.

Here we show that the converse problem, namely finding the least
common cover may have more than one solution. Let $G$ and $H$ be
disjoint connected graph. Let $G \smile H$ be a graph composed from
$G$ and $H$ adding an edge that connects some vertex of $G$ to some
other vertex of $H$. This operation depends on the choice of
vertices, bit we will indicate the choice simply by referring to the
figure. Let $G^*$ and $H^*$ be any two double covers of $G$ and $H$
respectively. Then there is a unique way to extend this to a double
cover of $G \smile H$ that we shall denote by $G^* \asymp H^*$.

Let us take the two familiar graphs $G(5,2)$ and $X$. Now form two
graphs $H_1 = G(5,2) \smile X \smile X$ and $H_2 = G(5,2) \smile
G(5,2) \smile X$; see Figure~\ref{fig_2}. Let $G_0 = G(10,3) \asymp
G(10,3) \asymp G(10,3), G_1 = G(10,3) \asymp G(10,3) \asymp 2X$ and
$G_2 = 2G(5,2) \asymp G(10,3) \asymp G(10,3).$

We claim that $G_0, G_1,$ and $G_2$ cover $H_1$ and $H_2.$ Using the
computer system {\sc Vega} (see \cite{Vega}) we checked that $G_0$
and $G_1$ are nonisomorphic covers of $H_1$ and $H_2$.

We may conclude by stating our finding in a more formal way.
\begin{Theorem}
There exist connected graphs $H_1$ and $H_2$ with a common universal
cover such that their minimal common cover is not unique.
\end{Theorem}

{\bf Acknowledgment.} The ideas of this work were conceived after
the Ledersprung Colloquium, which accompanies the traditional
Ledersprung, an annual student initiation event at the
Montanuniversit\"{a}t Leoben. The research was supported in part by
a grant P1-0294 from Ministrstvo za \v solstvo, znanost in \v sport
Republike Slovenije.

\end{document}